\documentstyle[12pt]{article}

  \newcommand{\const}{\rm const}

  \newcommand{\vraisup}{\rm vraisup}
  \newcommand{\Dom}{\rm  Dom}

\begin{document}

   \begin{center}

{\bf Tail estimates for random variables from interrelation between }\\

\vspace{4mm}

{\bf corresponding  moments inequalities. }\\

\vspace{5mm}

{\bf M.R.Formica, \ E.Ostrovsky, \ L.Sirota. } \\

 \end{center}

 \vspace{5mm}

\ Universit\`{a} degli Studi di Napoli Parthenope, via Generale Parisi 13, Palazzo Pacanowsky, 80132,
Napoli, Italy. \\
e-mail: \ mara.formica@uniparthenope.it \\

 \ Department of Mathematics and Statistics, Bar-Ilan University,\\
59200, Ramat Gan, Israel. \\
e-mail: \ eugostrovsky@list.ru\\

 \ Department of Mathematics and Statistics, Bar-Ilan University,\\
59200, Ramat Gan, Israel. \\
e-mail: \ sirota3@bezeqint.net \\

\vspace{5mm}

 \begin{center}

  {\bf Abstract.}

 \end{center}

\vspace{4mm}

 \hspace{3mm}  We derive the tail inequalities between two random variables starting from inequalities
 between its moment, or more generally between its Lebesgue - Riesz norms, which holds true on certain
 sets of parameters. \par
  \ We consider some applications into the martingale theory.\par

\vspace{5mm}

\hspace{3mm} {\it Key words and phrases.} \ Probability space, tail of distribution, random variable (r.v.), moment,
Young - Fenchel, or Legendre,  transform; generating function, natural function,
Lebesgue - Riesz and Grand Lebesgue Spaces and norms, martingale, Burkholder and Burkholder - Davis - Gundy inequalities,
estimations,  Lyapunov's inequality, examples.

\vspace{5mm}

\section{Introduction. Notations. Statement of problem.}

\vspace{5mm}

 \hspace{3mm} Let $ \ (\Omega = \{\omega\}, \ \cal{B}, \ {\bf P} )  \ $  be certain non - trivial  probability space with
an expectation $ \ {\bf E} \ $ and the classical Lebesgue - Riesz norm defined on the appropriate random variable (r.v.)
(measurable function)  $ \xi: \ \Omega \to R \ $

$$
|\xi|_p \stackrel{def}{=} \left[ \ {\bf E} |\xi|^p \ \right]^{1/p}, \ 1 \le p \le \infty,
$$
where as ordinary

$$
|\xi|_{\infty} := \vraisup_{\omega \in \Omega} |\xi(\omega)|.
$$
 \ Correspondingly, $ \ L_p = L_p(\Omega, {\bf P}) = \{ \ \xi: \ |\xi|_p < \infty. \ \} \ $ \par

\vspace{4mm}

 \ {\bf  We will start from the inequality of the form }

 \vspace{3mm}

\begin{equation} \label{key relation}
|\xi|_p \le g(p, r, |\eta|_r), \ (p,r) \in D.
\end{equation}

\vspace{3mm}
{\bf and intent to deduce the exact tail estimate for the first r.v.  } $ \ \xi \ $ {\bf via ones of the second one} $ \ \eta. \ $  \par

\vspace{4mm}

 \ Here $ \ D \ $ is certain  {\it non - empty} domain inside the quarter plane

$$
D \subset [1,\infty) \otimes [1,\infty),
$$
and $ \  g(p,r,z), \ p,r \ge 1, \ z \ge 0 \ $ is fixed numerical valued  deterministic non - negative
measurable function. \par

\vspace{4mm}

 \ A famous example: let $ \ (\epsilon_k, F_k), k = 1,2, \ldots,n \ $ be a centered martingale:

 $$
 {\bf E} \epsilon_m/F_k = \epsilon_k, \ k \le m.
 $$
 \ Set $ \ M^* = \max_k \epsilon_k, \  k \le n; \ $ then for  the values $ \ p > 1 \ $

\begin{equation} \label{Doob ineq}
|M^*|_p \le \frac{p}{p-1} \cdot |\epsilon_n|_p, \ -
\end{equation}
Doob's inequality. \par

\ This case was investigated in a previous article in this direction \cite{Ostrovsly 1}. Many other inequalities
of the form (\ref{key relation}) arises in the martingale theory: Burkholder - Davis - Gundy (BDG) inequality at so one, see e.g.  \cite{Burkholder 1},
\cite{Burkholder 2}, \cite{Burkholder Davis Gundy},  \cite{Hall Heyde}, \cite{Neveu}, \cite{Pena},  \cite{Revuz Yor} and so one. \par

\vspace{4mm}

 \ Another notations.

$$
R(p) = \{r: \ (p,r) \in D \}, \hspace{3mm} U = \{ \ p, \ \exists  r \ge 1, \ \Rightarrow (p,r) \in D \ \}.
$$

 \ Further, let $ \ p_0 \ $ arbitrary number from the set $ \ U. \ $ Define the following important function generated by the random
 variable $ \ \eta \ $  and by the values $  \ p \in U:  \ $

\begin{equation} \label{import fun upp}
  \psi(p) = \psi_{p_0}[\eta](p) \stackrel{def}{=} \inf_{r \in R(p)} g(p, r, |\eta|_r), \ p > p_0;
\end{equation}

\vspace{3mm}

\begin{equation} \label{import fun low}
\psi(p) = \psi_{p_0}[\eta](p) \stackrel{def}{=} \inf_{r \in R(p)} g(p_0, r, |\eta|_r), \ p \le p_0.
\end{equation}

 \vspace{4mm}

 \hspace{3mm}  We  have  from (\ref{key relation}) and from the Lyapunov's inequality again for the variables $ \ p \in U \ $

\vspace{4mm}

\ {\bf Proposition 1.1.}

\begin{equation} \label{Lp ineq}
|\xi|_p \le \psi(p) =  \psi_{p_0}[\eta](p), \ p \in U.
\end{equation}

\vspace{5mm}

\section{Grand Lebesgue Spaces estimates.}

\vspace{5mm}

 \hspace{3mm} The relation (\ref{Lp ineq}) may be rewritten  (and used) by means of the theory
 of the so - called Grand Lebesgue Spaces (GLS).  These spaces are investigated  in many works, see e.g.
\cite{Buldygin}, \cite{Capone1}, \cite{Capone2}, \cite{Ermakov},
\cite{Fiorenza-Formica-Gogatishvili-DEA2018},\cite{fioforgogakoparakoNAtoappear}, \cite{fioformicarakodie2017},
\cite{formicagiovamjom2015}, \cite{Formica Ostrovsky Sirota weak dep}, \cite{Iwaniec}, \cite{Kozachenko 1},
\cite{Kozachenko 2}, \cite{Liflyand}, \cite{Ostrovsky 0}, \cite{Ostrovsky 3}, \cite{Ostrovsky 4}. \par

\vspace{4mm}

 \hspace{3mm} We recall briefly the definition and some important properties of these GLS spaces, {\it adapted for this report.}
 Let $ \ \kappa = \kappa(p), \ p \in U \ $ be numerical valued positive measurable function, which may be named as {\it generating function}
 for introduced space.  Denote as ordinary by $ \ G\kappa \ $
 the set of all the numerical valued random variables $ \ \zeta \ $  having a finite norm

\vspace{3mm}

\begin{equation} \label{def Gpsi}
||\zeta||G\kappa \stackrel{def}{=} \sup_{p \in U}  \left\{ \ \frac{|\zeta|_p}{\kappa(p)} \ \right\}.
\end{equation}

 \ For instance, the {\it generating } function $ \ \kappa(\cdot) \ $ for this space may be introduced by the {\it natural} way:

$$
\kappa_0(p) = \kappa_0[\zeta](p) := |\zeta|_p,
$$
if of course it is finite for certain interval $ \ p \in (a,b), \ 1 \le a < b \le \infty. \ $ \par

\vspace{3mm}

 \ These  spaces are Banach functional complete and rearrangement invariant. They are closely related with the tail behavior
 of the r.v. $ \ \zeta: \ $

$$
T[\zeta](t) \stackrel{def}{=} {\bf P} (|\zeta| > t), \ t \ge 0.
$$

 \ Namely, introduce the auxiliary function

$$
S(x) = S[\psi](x) \stackrel{def}{=} \exp \left( \ - h^*[\psi](\ln x)  \ \right), \ x \ge e,
$$
where $ \ h(p) = h[\psi](p) = p \  \ln \ \psi(p),  \ $ and

\begin{equation} \label{Young Fenchel}
h^*(y) \stackrel{def}{=} \sup_{p \in \Dom[\psi]} (p y - h(p)) \ -
\end{equation}
is the famous Young - Fenchel, or Legendre,  transform for the function  $ \ h = h(p). \ $ \par

 \vspace{3mm}

  \ It is known \  \cite{Kozachenko 1},  \cite{Kozachenko 1},  \cite{Ostrovsky 0}, chapters 1,2  that if $ \ \zeta(\cdot) \in G\psi \ $ and
  $ \  || \zeta || = 1,   \ $ then

\begin{equation} \label{tail estim}
T[\zeta](x) \le S[\psi](x);
\end{equation}
and moreover the conversely proposition holds true. In detail, if the r.v. $ \ \zeta \ $ satisfies the estimate (\ref{tail estim}), then
under simple appropriate natural conditions on the function $ \ \psi  \ \Rightarrow \  \zeta \in G\psi: \ $

\begin{equation} \label{inverse}
\exists \ C = C[\psi] \ \in (0,\infty) \ \Rightarrow \ ||\zeta||G\psi  \le C(\psi).
\end{equation}

\vspace{3mm}

 \ For instance, the estimate  for some r.v. $ \ \zeta \ $  of the form

$$
\sup_{p \ge 1} \left\{ \ \frac{|\zeta|_p}{p^{1/m}}  \ \right\} < \infty,
$$
where $ \ m = \const > 0 \ $ is completely equivalent to the following tail estimate

$$
\exists \ c(m) \in (0,\infty) \ \Rightarrow  T[\zeta](t) \le \exp \left(-c(m) \ t^m \ \right), \ t \ge 0.
$$

\vspace{4mm}

 \ Consider in addition some comparison assertions. Let $ \ \nu_1 = \nu_1(p), \ \nu_2 = \nu_2(p), \ p \in U  \ $ be two functions
from the set $ \ G\Psi. \ $ Suppose

$$
\exists C \in (0,\infty) \ \Rightarrow  \forall p \in U \ \Rightarrow \nu_1(p) \le C \nu_2(p);
$$
then evidently

$$
||\zeta||G\nu_2 \le C \ ||\zeta||G\nu_1.
$$

\vspace{4mm}

 \ Further, suppose that for some r.v. $ \ \zeta, \ $ certain $ \ \psi \ - $  function $ \ \nu = \nu(p), \ p \in U \ $ and some
 fixed positive finite constant  $ \ C \in (0,\infty) \ $ there holds
$ \ ||\zeta||G\nu \le C, \  $  or equally

$$
|\zeta|_p \le C \ \nu(p), \ p \in U.
$$

 \ Then

\begin{equation} \label{C nu estim}
T[\zeta](t) \le \exp \left\{ \  - h^*[\nu] (\ln(t/C)) \ \right\},  \ t \ge C \ e.
\end{equation}

\vspace{3mm}

 \ To summarize: \par

\vspace{4mm}

  {\bf Proposition 2.1.} It follows under our notations and assumptions on the basis of the estimation  (\ref{Lp ineq})

\vspace{3mm}

\begin{equation} \label{main G psi}
T[\xi](t) \le \exp \left( \ - h^*[\psi]( \ln t) \ \right), \ t \ge e,
\end{equation}
with correspondent (exponential) tail estimate. \par

\vspace{5mm}

\section{Linear and  multilinear cases.}

\vspace{5mm}

 \hspace{3mm} Let us consider in this section a particular case of the inequality of the form

\begin{equation} \label{linear case}
|\xi|_p \le v(p,r) \ |\eta|_r,  \ (p,r)  \in D.
\end{equation}
  a "linear" case. This possibility appears in particular in  the mentioned above works devoted to the theory of martingales.\par

\vspace{3mm}

 \ Assume that the r.v. $ \ \eta \ $ belongs to some Grand Lebesgue Space $ \ G\beta: \ $

$$
|\eta|_r \le \beta(r) \ ||\eta||G\beta, \ r \in (a,b).
$$

 \ Introduce the new $ \ \psi \ - $ function

$$
\tau(p) := \inf_r [v(p,r) \ \beta(r)], \ p \in U.
$$

 \ We have $ \ |\eta|_r \le \beta(r) \ ||\eta||G\beta, \ $

$$
|\xi|_p \le v(p,r) \ \beta(r) \ ||\eta||G\beta,
$$
following

$$
|\xi(p) |_p \le \tau(p) \ \cdot ||\eta||G\beta.
$$

\vspace{3mm}

 \ To summarize: \par

 \vspace{4mm}

 \ {\bf Proposition 3.1.}

\begin{equation} \label{lin case}
||\xi||G\tau \le ||\eta||G\beta,
\end{equation}

\ with correspondent tail estimate.\par

\vspace{4mm}

 \ {\bf Remark 3.1.} The essential non - improvability in general case of the last estimate (\ref{lin case})
holds true for example in the case of Doob's  maximal inequality, see e.g.  \cite{Ostrovsly 1}; more simple
example: \ $ \ \xi = \eta. \ $\par

\vspace{3mm}

\ {\bf Remark 3.2.} The cases when $ \ |\xi|_p \le v_1(p,r) \ |\eta|_r^{\alpha} \ $ or when

$$
|\xi|_p \le v(p,\vec{r}) \ \prod_k |\eta_k|_{r(k)}^{\alpha(k)}
$$
may be investigated completely alike, as well as the case of the mixed (anisotropic) Lebesgue - Riesz norms

$$
|\eta|_{\vec{r}} = || \   || \ ||\eta||_{r_1, X_1} \ ||_{r_2, X_2} \ \ldots  \ ||_{r_d, X_d},
$$
where

$$
\eta = \eta(x_1, x_2, \ \ldots, \ x_d), \ d = 2,3,\ldots; \ x_j \in X_j,
$$

$ \ (X_j, \mu_j), j = 2,3,\ldots,d \ $ are  measurable spaces equipped correspondingly with the measures $ \ \mu_j. \ $ \par

\vspace{5mm}

\section{Examples.}

\vspace{5mm}

 \ {\bf A.} Doob's inequality, see (\ref{Doob ineq}).  Suppose for certain martingale  $ \ (\epsilon_k, F_k), \  0 \le k \le n, \ \epsilon_0 = 0 \ $

\begin{equation} \label{Doob Gpsim}
\exists C_2,m  \in (0,\infty)  \ \Rightarrow   T[\xi_n](t) \le \exp(- C_2 \ t^m), \ t \ge 0.
\end{equation}

 \ Then

\begin{equation} \label{Doob max m}
\exists C_3 = C_3(m) \in (0,\infty)  \ \Rightarrow   T[M^*n](t) \le \exp(- C_3(m) \ t^m), \ t \ge 0,
\end{equation}
 and conversely statement is also true. \par

 \vspace{3mm}

 \ {\bf B.} Burkholder - Davis - Gundy inequality. Let  $ \ (\xi(t), F_t), \ t \ge 0 \ $  be continuous time separable right - continuous
 martingale  for which $ \ \xi(0) = \xi(0+) = 0, \ $ (cadlag). Denote as ordinary

$$
M_{\infty}^*\stackrel{def}{=} \sup_{t \ge 0} \xi(t),
$$
and by $ \ <M,M> \ $ its quadratic variation, if both this variables there exists (and are finite). The famous Burkholder - Davis - Gundy
upper, (as well as lower,) inequality has a form

$$
| \ M_{\infty}^* \ |_p \le \sqrt{e} \cdot | \  <M.M> \ |_{p/2}^{1/2}, \ p \ge 2.
$$

 \ One can conclude as before that if

$$
T[\sqrt{<M,M>}](t) \le \exp(- c_4 \ t^m),  \ t > 0, \ \exists \ c_4 \in (0,\infty),
$$
then alike

$$
T[M^*_{\infty}](t) \le \exp(- c_5 \ t^m),  \ t > 0, \ \exists \ c_5 = c_5(m) \in (0,\infty),
$$
and conversely statement is also true. \par

\vspace{6mm}

\vspace{0.5cm} \emph{Acknowledgement.} {\footnotesize The first
author has been partially supported by the Gruppo Nazionale per
l'Analisi Matematica, la Probabilit\`a e le loro Applicazioni
(GNAMPA) of the Istituto Nazionale di Alta Matematica (INdAM) and by
Universit\`a degli Studi di Napoli Parthenope through the project
\lq\lq sostegno alla Ricerca individuale\rq\rq .\par

\end{document}